\documentclass[12pt]{article}
\usepackage{geometry}
\usepackage{graphicx}
\usepackage{amsmath, amssymb, amsfonts}
\usepackage{algorithm, algpseudocode}
\usepackage{bm}
\usepackage{tikz}
\usepackage{stmaryrd}
\usepackage{appendix}
\usepackage{lineno}
\usepackage[affil-it]{authblk}
\usepackage[numbers, sort, comma, square]{natbib}

\algrenewcommand\algorithmicfor{}
\algrenewcommand\algorithmicend{{\bf{End do}}}

\providecommand{\keywords}[1]
{
  \small	
  \textbf{\textit{Keywords---}} #1
}

\title{Exploring physics-dynamics coupling using moist shallow water equations}

\author[1,*]{Nell Hartney}
\author[2]{Thomas M. Bendall}
\author[1]{Jemma Shipton}

\affil[1]{Department of Mathematics and Statistics, University of Exeter, Exeter, UK}
\affil[2]{Dynamics Research, Met Office, Exeter, UK}
\affil[*]{Corresponding author: Nell Hartney, nh491@exeter.ac.uk}

\begin{document}

\maketitle

\begin{abstract}
    One of the key choices for numerical models of geophysical fluids is how parametrisations of physical processes interact with the numerical methods that handle the resolved flow, known in the atmospheric community as the dynamical core. As both the dynamical core and parametrisations of physics processes continue to evolve and improve, the issue of physics-dynamics coupling - how these two different parts of the model interact - becomes ever more important. In this paper we use two variations of the moist shallow water equations to develop a simplified framework that can be used to investigate some of the questions associated with physics-dynamics coupling. The shallow water equations act as a simplified dynamical core that is computationally cheap but still retains pertinent features of the atmosphere, and the introduction of moisture means the addition to the model of a physical parametrisation. This study uses `split-physics' moist thermal shallow water equations which couple moisture to the shallow water equations via a parametrisation, and also develops a new `integrated-physics' formulation of the moist thermal shallow water equations which re-formulates the model so that all the moist processes are captured in the dynamics. This integrated-physics model thus requires no physics-dynamics coupling and acts as a ground-truth to compare coupling strategies in the split-physics formulation against. We use both models to examine the effect of varying the approach to how physics is coupled to the dynamics in the semi-implicit quasi-Newton timestepping scheme. The results demonstrate the usefulness of the integrated-physics moist shallow water equations and provide insights into how best to deal with physics in a model's timestep.
\end{abstract}


\section{Introduction}
A significant challenge in applying the techniques of fluid dynamics to the atmosphere and ocean is the range of spatial and temporal scales that atmospheric and oceanic processes span. For example, processes such as turbulence typically have scales of a few seconds in time and metres in space, while processes like planetary waves can stretch for thousands of kilometres and last for several days. This means a weather or climate forecast must deal with a number of physical processes at very different space and time scales simultaneously - a significant challenge in developing solution strategies for these equations. The result of this range of scales is a collection of processes that cannot be resolved by the numerical model even though they may be critical to producing an accurate forecast, interacting with and influencing the large-scale dynamics as they do. For this reason these unresolved processes are usually parametrised - incorporated with an alternate mathematical model - to allow their effect on the overall flow to be included correctly. These parametrised processes are known as the \textit{physics} and the fluid dynamics that are discretised by the dynamical core the \textit{dynamics}. As well as processes that are unresolved on the scale of the dynamical core, the physics includes non-fluid dynamical elements that also influence the weather and climate such as radiation. Producing a full coherent model requires coupling together the separate components of physics and dynamics in some manner, a task which has become known as \textit{physics-dynamics coupling}.

There are numerous challenges associated with physics-dynamics coupling, and, as described in \cite{PDCreview}, the implementation of the coupling can have a significant impact on the model. Improvements in individual model components can easily be lost if the errors introduced by the coupling dominate, and hence physics-dynamics coupling is one of the important factors that limits the overall accuracy of a full model. It is for this reason that the weather and climate modelling communities have, relatively recently, started to make focused efforts in addressing physics-dynamics coupling as a topic in its own right.

The complexity of the dynamical core in a weather model and the range of physics processes that must be included to produce an accurate forecast makes investigating physics-dynamics coupling questions in a full weather model very challenging. There is also the problem that it is difficult to evaluate the effectiveness of coupling strategies when there is no way to produce a `true' solution for comparison. This paper aims to address both of these issues through the use of a simplified model with a single physics process, and a re-formulation of the same model that requires no coupling and can act as a reference solution. 

The simplified model we have chosen is the moist shallow water model, which introduces moisture to the traditional shallow water equations. The shallow water equations in their usual form are widely used in the testing of new algorithms for weather and climate modelling. The primary reason for this is that while computationally cheap, the shallow water equations still represent pertinent features of the atmosphere such as Rossby waves and geostrophic adjustment, and often represent an optimal compromise between complexity and an accurate representation of dynamical features. Like other simplified models, they are based on systematic approximations of the full equations and elimination of certain spatial and temporal scales \cite{Zeitlinbook}. The shallow water equations exploit the \textit{shallowness} of the atmosphere and ocean; that is, the fact that the atmosphere and the oceans are both very thin with respect to the radius of the Earth. In their traditional form, however, the shallow water equations are a purely dynamical model and do not include any physics parametrisations. For our physics-dynamics coupling investigations, therefore, we will use an adaptation of the shallow water equations that includes moist physics. This gives us a model that is simple and computationally cheap but also includes a physics process, thus giving us the opportunity to investigate questions about time-stepping with physics and physics-dynamics coupling questions. The choice of moisture as a physics process is motivated by the fact that moist processes generally include a discontinuity, which is not only numerically challenging but is also a common feature of many physics schemes. Dealing with discontinuities here makes the approach relevant to other physics parametrisations that include this typical numerical feature. The second reason we have chosen to work with moisture as a physics processes is that we can exploit the idea of instantaneous saturation adjustment common in many moist parametrisations to create the `integrated-physics' formulation of the equations. This formulation acts as a control and allows us to quantify errors in different coupling approaches. 

Different variants of the moist shallow water equations have been proposed by different authors (see, for example, \cite{Gill1982, Bouchut2009, Lambaertsetal2011baroclinic, KLZ2020b, ZA2015, Yang2021, WC2014}). The most commonly-used moist shallow water models add equations for moist variables to the standard shallow water set, and then couple them to the pre-existing equations in various ways. A general framework for these type of moist shallow water equations was outlined in \cite{Hartney2025}, which encapsulated the moist convective shallow water equations, the moist thermal shallow water equations and the moist convective thermal shallow water equations, in addition to a new fourth formulation called the moist convective pseudo-thermal shallow water equations. Common to all these approaches recoverable from the general one is that fact that each model includes physics; that is, the equations have source terms which are dictated by the moist physics scheme. It is these source terms that necessitate a physics-dynamics coupling strategy, and mimic source terms from physics schemes in a full weather model.

The novelty of this paper lies in the development of a moist shallow water model with no physics source terms, which removes the need for a physics-dynamics coupling strategy. Removing any coupling means that this dynamics-only, integrated-physics model can act as a `ground-truth' to compare coupling approaches against. We use this reference solution in a moist gravity wave test case to investigate different approaches to physics in the semi-implicit quasi-Newton time-stepping scheme. This time-stepping scheme follows the approach of the Met Office's ENDGame \cite{Wood2014inherently} and GungHo \cite{GungHoCartesian} dynamical cores, and our models are discretised in space using the compatible finite element method, which will form the basis of the Met Office's next generation dynamical core. These choices are motivated by our desire for the investigation to be as relevant as possible to this modern operational numerical weather prediction model. Though this integrated-physics idea would not be possible in a full dynamical core with more sophisticated dynamics (and thus there is no full-model analogue for what we are doing), our belief is that the learnings around coupling in this simplified shallow water framework can be useful in more complex frameworks too.  

The structure of the article is as follows: in Section 2 we give an overview of the `split-physics' moist thermal shallow water model. The split-physics model acts as our coupled model that we will investigate coupling strategies in. Section 2 also introduces the `integrated-physics' formulation of the moist thermal shallow water equations, and shows how it is derived from the split-physics variant of the same model. This section thus provides us with the framework that we can use to compare various different coupling approaches to the reference solution. Both the compatible finite element discretisations of the models and the semi-implicit quasi-Newton time-stepping scheme are dealt with in Section 3, which demonstrates how our approach to solving the equations is highly relevant in a modern numerical weather prediction paradigm. Finally, the two test cases and results in Section 4 verify that both of the model formulations behave as expected and show the effect of varying how physics interacts with the semi-implicit quasi-Newton timestep.

\section{Equations}
\subsection{The shallow water equations}
The shallow water equations in their traditional form are a purely dynamics model, with no physics parametrisations. They are derived, as described in, for example \citet{Zeitlinbook}, by integrating the primitive equations between two material surfaces, assuming that the horizontal flow is independent of depth and applying the hydrostatic balance condition. Written here in advective form with rotation and without any dissipative effects, the shallow water equations are given by:
\begin{align}
    &\frac{\partial{\bm{u}}}{\partial{t}} + (\bm{u} \cdot \bm{\nabla}) \bm{u}  + f \bm{\hat{k}} \times \bm{u} + g \bm{\nabla} (D + B) = 0, \\
    &\frac{\partial{D}}{\partial{t}} + \bm{\nabla} \cdot (\bm{u}D) = 0,
\end{align}
where $\bm{u}$ is horizontal velocity, $D$ is the thickness or depth of the fluid layer, $g$ is acceleration due to gravity, $f$ is the Coriolis parameter, $\bm{\hat{k}}$ is the unit vector in the vertical and $B$ is the height of the bottom surface.

Introducing moisture into the shallow water system provides an example of a physics parametrisation, which can be coupled to the dynamics in various ways by choosing the impact that moist phase changes should have on the existing shallow water dynamics. One possibility is to model the thermal effect of moist phase changes using the moist thermal shallow water equations. These equations were proposed by \citet{ZA2015}, where they derived a moist shallow water system from the three-dimensional Boussinesq approximation. The thermal variable in their model is the potential temperature, which  was converted to buoyancy in \citet{Hartney2025}. It is this buoyancy formulation of the moist thermal shallow water equations that we use in this work.

\subsection{The split-physics moist thermal shallow water equations}
We will refer to this buoyancy formulation of the moist thermal shallow water equations as the \textit{split-physics} formulation of the model, because this formulation couples a moist physics scheme to the thermal shallow water equations using source terms, and thus has the physics `split' from the dynamics. The equations are given by:
\begin{align}
    &\frac{\partial \bm{u}}{\partial t}
    + (\bm{u} \cdot \bm{\nabla}) \bm{u} 
    + f \bm{\hat{k}} \times \bm{u} + b \bm{\nabla} (D + B) + \frac{D}{2} \bm{\nabla} b = 0, \label{eq:split_phys_u}\\
    &\frac{\partial D}{\partial t} + \bm{\nabla} \cdot ({\bm{u} D}) = 0, \label{eq:split_phys_D} \\
    &\frac{\partial b}{\partial t} + (\bm{u} \cdot \bm{\nabla}) b = -\beta_2 P, \label{eq:split_phys_b}\\
    &\frac{\partial q_v}{\partial t} + (\bm{u} \cdot \bm{\nabla}) q_v = -P, \label{eq:split_phys_qv} \\
    &\frac{\partial q_c}{\partial t} + (\bm{u} \cdot \bm{\nabla}) q_c = P, \label{eq:split_phys_qc}
\end{align}
where the same shallow water notation applies as above, $b$ is the buoyancy, $q_v$ is the water vapour and $q_c$ is the cloud water. $P$ is the source term from the physics scheme, described in Section \ref{section:moist_physics_scheme} below, that arise as moist species undergo phase changes. The moisture variables $q_v$ and $q_c$ represent mixing ratios, as in Zerroukat and Allen's moist thermal shallow water model \citep{ZA2015}. $\beta_2$ is $gL$ where $L \approxeq 10$ is the pseudo-latent heat from \citet{ZA2015}.

\subsection{The integrated-physics moist thermal shallow water equations: a ground truth}
The aim of the \textit{integrated-physics} formulation of the moist thermal shallow water equations is to re-formulate the split-physics version of the equations to move all the moist physics into the dynamics. This means a model with no physics source terms, which removes the need for a physics-dynamics coupling strategy. Removing any coupling means that this dynamics-only, integrated-physics model can act as a `ground-truth' to compare coupling approaches in the split-physics formulation against.

We begin with the split-physics moist thermal shallow water equations given by Equations \eqref{eq:split_phys_u} - \eqref{eq:split_phys_qc}. We define the equivalent buoyancy, $b_e$, as $b_e = b - \beta_2 q_v$. This is a conserved quantity that will replace the buoyancy in the split-physics moist thermal shallow water equations. We also define a total moisture variable, $q_t$, which is the sum of the water vapour and the cloud. The model assumes that condensation and evaporation are instantaneous processes so that $q_v$ is always equal to saturation, if there is enough moisture.
The advection equations for water vapour and cloud in the split-physics formulation (Equations \eqref{eq:split_phys_qv} and \eqref{eq:split_phys_qc}) are replaced in the integrated-physics formulation by an advection equation for the (conserved) total moisture $q_t$. Then the integrated-physics formulation is given by:
\begin{align}
    &\frac{\partial \bm{u}}{\partial t}
    + (\bm{u} \cdot \bm{\nabla}) \bm{u}
    + f \bm{\hat{k}} \times \bm{u} + (b_e + \beta_2 q_v) \bm{\nabla} (D + B) + \frac{D}{2} \bm{\nabla} (b_e + \beta_2 q_v) = 0, \label{eq:int_phys_u} \\
    &\frac{\partial D}{\partial t} + \bm{\nabla} \cdot ({\bm{u} D}) = 0, \label{eq:int_phys_D} \\
    &\frac{\partial b_e}{\partial t} + (\bm{u} \cdot \bm{\nabla}) b_e = 0, \label{eq:int_phys_b} \\
    &\frac{\partial q_t}{\partial t} + (\bm{u} \cdot \bm{\nabla}) q_t = 0. \label{eq:int_phys_qt}
\end{align}

We now have $q_v$ appearing in the pressure gradient-like term on the right-hand side of the $\bm{u}$ equation (Equation \eqref{eq:int_phys_u}). Since we transport only $q_t$ and not $q_v$ itself, we must diagnose $q_v$ from $q_t$ to adjust the $\bm{u}$ equation at every timestep. To do this we use a saturation function, $q_{sat}$, to partition $q_t$ into $q_v$ and $q_c$, where all moisture below the saturation function is water vapour and all moisture above it is cloud: 
\begin{equation} \label{eq:md_sat_function}
    q_{sat} = \frac{q_0 H}{D + B} e^{\nu(1-b_e/g)},
\end{equation}
where $q_0$ and $\nu$ are scaling factors (discussed in the results section), $H$ is the background depth of the fluid layer and $g$ is acceleration due to gravity. We have chosen this as a semi-realistic saturation function which is a function of $b_e$ but not $q_v$, which simplifies the model. We use the saturation function at each timestep to diagnose $q_v$, via:
\begin{equation}
    q_v = \begin{cases}
        q_{sat}, & \text{if } q_t \geq q_{sat} \\
        q_t, & \text{otherwise}
    \end{cases}
\end{equation}
This $q_v$ is then used in the computation of the pressure gradient-like term on the right-hand side of the $\bm{u}$ equation (Equation \eqref{eq:int_phys_u}).

\subsection{Moist physics scheme} \label{section:moist_physics_scheme}
For this work the split-physics equations are coupled to a moist physics scheme with reversible conversions between two species: water vapour $q_v$ and cloud water $q_c$. The scheme is based on the three-state moist physics scheme described in \citep{Hartney2025}, though here we exclude rain to model just a reversible phase change between water vapour and cloud. Conversions between the two are dictated by a prescribed saturation function $q_{sat}$:
\begin{equation}
    q_{sat} = \frac{q_0 H}{D + B} e^{20(1-b/g - \beta_2 q_v/g)}, \label{eq:sat_func}
\end{equation}
where $q_0$ is a scaling factor discussed below in the results section, $H$ is the background depth and $B$ is the topography. Whenever the local value of the vapour $q_v$ exceeds the saturation function then a portion of the excess $(q_v - q_{sat})$ will condense into cloud over a timescale $\tau_v$, causing a local heating of the atmosphere. In this work we choose the condensation timescale to be the same as the timestep such that $\gamma_v = \Delta t$, so that moist processes are considered instantaneous. This conversion can be written as:
\begin{align}
    q_v > q_{sat} : P &= \gamma_v \frac{q_v - q_{sat}}{\Delta t}, \label{eq:evap_P}
\end{align}
where $\gamma_v$ is a conversion rate that depends on the depth and buoyancy and is discussed below.

In a similar way, in a sub-saturated atmosphere where cloud is present, cloud will evaporate with a local cooling effect:
\begin{align}
    q_v < q_{sat} : - P &= \text{min}\left [\frac{q_c}{\Delta t}, \gamma_v \frac{q_{sat} - q_v}{\Delta t} \right ], \label{eq:condens_P}
\end{align}
where the minimum function is used to ensure that only the available cloud can be evaporated.

These conversions adjust the buoyancy, water vapour and cloud equations via:
\begin{align}
    &b \rightarrow b - \beta_2 P, \label{eq:b_increment} \\
    &q_v \rightarrow q_v - P, \label{eq:v_increment} \\ 
    &q_c \rightarrow q_c + P.
\end{align}

The factor $\gamma_v$ ensures that only a fraction of the cloud or vapour is converted each time a phase change happens, which prevents a two-step oscillation arising where a phase change causes a change in temperature which subsequently returns the atmosphere to the starting state. This factor is derived in \citep{Hartney2025} and is given by:  
\begin{equation} \label{eq:gamma_v}
    \gamma_v = \frac{1}{1 + q_{sat}(D,b)\left[ \frac{20 \beta_2}{g}\right]}.
\end{equation}  

Putting these two conversion processes together (Equations \eqref{eq:evap_P} and \eqref{eq:condens_P}) the source term $P$ can be concisely written as:
\begin{equation}
    P = \text{max}\left(\gamma_v\frac{q_v-q_{sat}}{\Delta t}, -\frac{q_c}{\Delta t}\right).
\end{equation}

\section{Discretisation}
\subsection{Compatible finite element discretisation}
Both the split-physics and the integrated-physics moist thermal shallow water models are discretised in space using the compatible finite element method. One of the attractions of compatible finite element methods is that they offer many of the desirable properties analogous to the C-grid for finite difference methods, but are suitable on non-orthogonal grids (see, for example, \citet{Cotter2023} for a discussion). It is for this reason that the Met Office will use compatible finite element methods as the basis for its new dynamical core, facilitating as it does the move away from the latitude-longitude grid to the cubed-sphere. The cubed-sphere grid is quasi-uniform over the sphere and so addresses the issues with parallel scalability that arise with the latitude-longitude grid as supercomputer architectures migrate towards massively parallel systems.

The essence of the finite element method is to approximate fields as the sum of a finite series of functions defined on a mesh of cells, and write the equations in weak form by multiplying by test functions and integrating over the domain. Compatible finite element methods use function spaces based on a discrete de Rham complex, such that the differential operators $\bm{\nabla}$, $\bm{\nabla} \times $ and $\bm{\nabla} \cdot$ map surjectively onto the kernel of the next operator in the sequence \citep{Cotter2023}. Cotter and Shipton first demonstrated in \citep{CotterShipton2012} that it is this compatible formulation that is key to maintaining the desirable properties of the C-grid, such as the ability to support steady geostrophic modes, the avoidance of spurious mode branches, and an accurate representation of the dispersion relation for Rossby and inertia-gravity waves. Here we use an extension to the discretisation given first in \citet{CotterShipton2012} to include finite element spaces for buoyancy and moisture variables, described in \citet{Hartney2025}. The discretisation is based on a two-dimensional discrete de Rham complex $(\mathbb{V}^1, \mathbb{V}^2)$, where $\mathbb{V}^1 = BDM_2$ (Brezzi-Douglas-Marini elements on triangles) and $\mathbb{V}^2 = DG_1$. We use the $\mathbb{V}^1$ space for the $\bm{u}$ variable and the $\mathbb{V}^2$ for the buoyancy and moisture variables. This means that the $\bm{u}$ belongs in a space of quadratic vector-valued functions whose normal components are continuous across cell edges, and the other prognostic fields in a space of linear functions with no inter-element continuity constraints - the next-to-lowest order spaces for this finite element family.

\subsubsection{Split-physics equations}
We write the weak form of the $\bm{u}$ equation in the split-physics formulation by multiplying Equation \eqref{eq:split_phys_u} by a test function ${\bm{\psi}}$ from the $\mathbb{V}^1$ space and integrating. After integration by parts we have:
\begin{equation}
\begin{split}
    \int_{\Omega} {\bm{\psi}} \cdot \frac{\partial {\bm{u}}}{\partial t} \text{d}x
    +\int_{\Omega} {\bm{\psi}} \cdot (\bm{u} \cdot \bm{\nabla}) \bm{u} \text{d}x
    + \int_{\Omega} {\bm{\psi}} \cdot (f \bm{\hat{k}} \times \bm{u}) \text{d}x
    - \int_{\Omega} (B + D) \bm{\nabla} \cdot ({\bm{\psi}}b) \text{d}x
    \\
    + \int_{\Gamma} \langle B + D\rangle \llbracket b{\bm{\psi}} \rrbracket \text{d}S
    - \frac{1}{2} \int_{\Omega} b \bm{\nabla} \cdot (D {\bm{\psi}}) \text{d}x
    + \frac{1}{2} \int_{\Gamma} \langle b\rangle \llbracket D{\bm{\psi}} \rrbracket \text{d}S = 0, \quad \forall \bm{\psi} \in \mathbb{V}^1.
\end{split}
\end{equation}
Here $\Gamma$ represents the set of cell edges and $\text{d}S$ is the measure for integrating over all cell edges. The definition of discontinuous fields on facets is taken as the product of the `jump' across the facet (denoted by $\llbracket \psi \rrbracket = \psi^+ - \psi^-$, where each edge has sides arbitrarily labelled ``+" and ``-") and the average (denoted by $\langle \psi\rangle = (\psi^+ + \psi^-)/2$).

We write the weak form of the depth equation (Equation \eqref{eq:split_phys_D}) by multiplying by a test function $\phi$ from the $\mathbb{V}^2$ space and integrating, using integration by parts on the second term. Since $D$ is being transported, we define $D$ on cell edges as the upwind value with respect to the velocity, denoted by $\Tilde{D}$. This gives the weak form:
\begin{equation}
    \int_{\Omega} \phi \frac{\partial D}{\partial t} \text{d}x - \int_{\Omega} D {\bm{u}} \cdot  \bm{\nabla} \phi  \text{d}x + \int_{\Gamma} \Tilde{D} \llbracket \phi \rrbracket {\bm{u}} \cdot \bm{n} \text{d}S = 0, \quad \forall \phi \in \mathbb{V}^2.
\end{equation}
In a similar way, the weak form of the buoyancy equation (Equation \eqref{eq:split_phys_b}) is:
\begin{equation}
    \int_{\Omega} \lambda \frac{\partial b}{\partial t} \text{d}x - \int_{\Omega} b \bm{\nabla} \cdot (\lambda {\bm{u}}) \text{d}x + \int_{\Gamma} \Tilde{b} \llbracket \lambda {\bm{u}}\rrbracket \text{d}S = \int_{\Omega} \lambda \beta_2 P \text{d}x, \quad \forall \lambda \in \mathbb{V}^2,
\end{equation}
and the weak form of the water vapour and cloud equations (Equations \eqref{eq:split_phys_qv} and \eqref{eq:split_phys_qc} respectively) is:
\begin{equation}
    \int_{\Omega} \tau \frac{\partial q^{(k)}}{\partial t} \text{d}x - \int_{\Omega} q^{(k)} \bm{\nabla} \cdot (\tau {\bm{u}}) \text{d}x + \int_{\Gamma} \Tilde{q}^{(k)} \llbracket \tau {\bm{u}} \rrbracket \text{d}S = \pm \int_\Omega \tau P \text{d}x, \quad \forall \tau \in \mathbb{V}^2,
\end{equation}
where $q^{(k)}$ represents each moist species, and $P$ on the right-hand sign is positive for water vapour and negative for cloud. 

\subsubsection{Integrated-physics equations}
As with the split-physics formulation of the equations, the weak form of the $\bm{u}$ equation in the integrated-physics formulation is derived by multiplying by a test function ${\bm{\psi}}$ from the $\mathbb{V}^1$ space. This gives the weak form for Equation \eqref{eq:int_phys_u} (after integration by parts) as:
\begin{equation}
\begin{split}
    \int_{\Omega} {\bm{\psi}} \cdot \frac{\partial {\bm{u}}}{\partial t} \text{d}x
    +\int_{\Omega} {\bm{\psi}} \cdot (\bm{u} \cdot \bm{\nabla}) \bm{u} \text{d}x
    + \int_{\Omega} {\bm{\psi}} \cdot (f \bm{\hat{k}} \times \bm{u}) \text{d}x
    - \int_{\Omega} (B + D) \bm{\nabla} \cdot ({\bm{\psi}}(b_e + \beta_2 q_v)) \text{d}x
    \\
    + \int_{\Gamma} \langle B + D\rangle \llbracket {\bm{\psi}} (b_e + \beta_2 q_v) \rrbracket \text{d}S
    - \frac{1}{2} \int_{\Omega} (b_e + \beta_2 q_v) \bm{\nabla} \cdot (D {\bm{\psi}}) \text{d}x
    \\
    + \frac{1}{2} \int_{\Gamma} \langle b_e + \beta_2 q_v \rangle \llbracket D{\bm{\psi}} \rrbracket
    \text{d}S 
    = 0 , \quad \forall \bm{\psi} \in \mathbb{V}^1.
\end{split}
\end{equation}
Again here $\Gamma$ represents the set of cell edges and $\text{d}S$ is the measure for integrating over all cell edges. As in the split-physics weak form, the definition of discontinuous fields on facets is taken as the product of the `jump' across the facet and the average.

The weak form of the depth equation (Equation \eqref{eq:int_phys_D} is derived by multiplying by a test function $\phi$ from the $\mathbb{V}^2$ space and integrating, using integration by parts on the second term. As in the split-physics formulation, we define $D$ on cell edges as the upwind value with respect to the velocity because it is being transported. This gives the weak form:
\begin{equation}
    \int_{\Omega} \phi \frac{\partial D}{\partial t} \text{d}x - \int_{\Omega} D {\bm{u}} \cdot \bm{\nabla} \phi  \text{d}x + \int_{\Gamma} \Tilde{D} \llbracket \phi \rrbracket {\bm{u}} \cdot \bm{n} \text{d}S = 0, \quad \forall \phi \in \mathbb{V}^2.
\end{equation}
In a similar way, the weak form of the equivalent buoyancy equation (Equation \eqref{eq:int_phys_b}) is:
\begin{equation}
    \int_{\Omega} \lambda \frac{\partial b_e}{\partial t} \text{d}x - \int_{\Omega} b_e \bm{\nabla} \cdot (\lambda {\bm{u}}) \text{d}x + \int_{\Gamma} \Tilde{b_e} \llbracket \lambda {\bm{u}}\rrbracket \text{d}S = 0, \quad \forall \lambda \in \mathbb{V}^2,
\end{equation}
and the moisture equation (Equation \eqref{eq:int_phys_qt}) is:
\begin{equation}
    \int_{\Omega} \tau \frac{\partial q_t}{\partial t} \text{d}x - \int_{\Omega} q_t \bm{\nabla} \cdot (\tau {\bm{u}}) \text{d}x + \int_{\Gamma} \Tilde{q_t} \llbracket \tau {\bm{u}} \rrbracket \text{d}S = 0, \quad \forall \tau \in \mathbb{V}^2.
\end{equation}

\subsection{Semi-implicit quasi-Newton timestepper}
\label{section:SIQN_scheme}
The equations are discretised in time using a semi-implicit quasi-Newton (SIQN) time-stepping scheme. This scheme follows the time stepper used by the Met Office's models, as described in \citet{Wood2014inherently} and \citet{GungHoCartesian}, and is also the time stepper used in \citet{Bendall2019recovered} and \citet{Bendall2020compatible}. The scheme is outlined in pseudocode in Algorithm \ref{algorithm:SIQN_alg}, and is also described in detail in \citet{Hartney2025}.

\begin{algorithm}
	\caption{Semi-implicit quasi-Newton time-stepping scheme}
    \label{algorithm:SIQN_alg}
	\begin{algorithmic}[1]
        \State Set: $\chi^{n+1} = \chi^n$
        \State Explicit forcing: $\chi_{f,e} = \chi^n - \frac{1}{2} \Delta t \mathcal{F}(\chi^n)$
		\State Outer loop: \textbf{do}
            \State \hspace{\algorithmicindent} Update advecting velocity: $\bar{\mathbf{u}} = \frac{1}{2}(\mathbf{u}^{n+1} + \mathbf{u}^n)$
            \State \hspace{\algorithmicindent} Explicit transport: $\chi_T = \mathcal{A}_{\bar{\mathbf{u}}}(\chi_{f,e})$
            \State \hspace{\algorithmicindent} Inner loop: \textbf{do}
                    \State \hspace{\algorithmicindent}\hspace{\algorithmicindent} Implicit forcing: $\chi_{f,i} = \chi_T - \frac{1}{2} \Delta t \mathcal{F}(\chi^{n+1})$
                    \State \hspace{\algorithmicindent}\hspace{\algorithmicindent} Write residual: $\chi_r = \chi_{f,i} - \chi^{n+1}$
                    \State \hspace{\algorithmicindent}\hspace{\algorithmicindent} Solve: $\mathcal{S}(\delta \chi) = \chi_r$ for $\delta \chi$
                    \State \hspace{\algorithmicindent}\hspace{\algorithmicindent} Increment: $\chi^{n+1} = \chi^{n+1} + \delta \chi$
            \State \hspace{\algorithmicindent} \textbf{End do}
    \State \textbf{End do}
    \State Explicit physics: $\chi^{n+1} = \mathcal{P}(\chi^{n+1})$
    \State Advance time step: $\chi^n = \chi^{n+1}$
	\end{algorithmic} 
\end{algorithm}

The stepper begins with a half timestep of explicit forcing which takes the state at the $n$th time level, $\chi^n$, and returns $\chi_{f, e}$. The forcing terms included both in this explicit forcing step and later in the implicit forcing step are the Coriolis and pressure gradient terms, which in Gusto are all written on the left-hand side. This results in a minus sign in the forcing computation in the SIQN algorithm.

The next step is an outer iterative loop, which begins with computation of the advecting velocity $\bar{\bm{u}}$ as the average of the velocity at the current time level $\bm{u}^{n}$ and the best guess of the velocity at the next time level $\bm{u}^{n+1}$. This advecting velocity is then used to explicitly advect the prognostic fields, $\mathcal{A}_{\bar{\bm{u}}}(\chi_{f, e})$, returning $\chi_T$. We use the advective form of the wind equation and transport all fields using an upwind discontinuous Galerkin scheme with a three-step Runge-Kutta timestepping procedure (SSPRK3), as outlined in, for example \citep{Bendall2020compatible} and \citep{Bendall2019recovered}. After that the algorithm enters the inner loop which begins with a half timestep of implicit forcing. We define a residual $\chi_{r}$ as the difference between this newly-forced state and the best guess for the state at the next timestep. We then solve a linear system iteratively, with the goal of reducing the residual to zero. Rather than computing the Jacobian at every iteration $k$ we use a quasi-Newton approach, which approximates the Jacobian as a simple linear system. Then the linear solve step seeks solutions $\delta \chi$ to the linear problem:
\begin{equation} \label{eq:linear_solve}
    \mathcal{S} \delta \chi = - \chi_{r},
\end{equation}
where $\mathcal S$ is the linear operator, $\delta \chi$ is the increment $\chi_{k+1}^{n+1} - \chi_{k}^{n+1}$ and $\chi_r$ is the residual.

The choice of this linear operator $\mathcal{S}$ is inspired by the linearisation of the residual $\chi_r$ about some reference state $\chi_{ref}$, following \citet{Wood2014inherently}. In this work we use two different approaches to this linear problem: firstly using a dry linear solver and secondly using a moist linear solver. The dry linear solver is derived by linearising the thermal shallow water equations about a reference state. It excludes moist variables (following the approach taken in the Met Office's ENDGame \citep{Wood2014inherently} and GungHo \citep{GungHoCartesian} models and the compressible Euler model from \citet{Bendall2020compatible}). In contrast, the moist linear solver includes moisture as prognostic variables in the linear solve. The dry linear solver for the split-physics version of the equations is outlined in detail in \citet{Hartney2025}, so we describe only the integrated-physics version of the dry solver and the moist solver here.

\subsubsection{Integrated-physics formulation of the dry linear solver}
We formulate the integrated-physics version of the dry linear solver by linearising the integrated physics equations (excluding topography) about a reference state $(\bm{0}, \bar{D}, \bar{b_e}, \bar{q_v})$. Written in weak form and after integrating by parts, this gives the system to solve (Equation \eqref{eq:linear_solve}) as:

\begin{align}
\begin{split}
    \int_\Omega {\bm{\psi}} \cdot \delta \bm{u} \text{d}x
    + \alpha \Delta t \int_{\Omega} {\bm{\psi}} \cdot (f \bm{\hat{k}} \times \delta \bm{u}) \text{d}x
    - \alpha \Delta t \int_{\Omega} \bm{\nabla} \cdot \left(\bm{\psi}(\bar{b_e}
    + \beta_2 \bar{q_v})\right) \delta D \text{d}x
    \\
    + \alpha \Delta t \int_{\Gamma} \llbracket (\bar{b_e} + \beta_2 \bar{q_v}) \bm{\psi} \rrbracket \langle \delta D\rangle \text{dS}
    - \frac{\alpha \Delta t}{2}\int_{\Omega} \delta b_e \bm{\nabla} \cdot (\bar{D} \bm{\psi}) \text{d}x
    - \frac{\alpha \Delta t}{2} \int_{\Omega} (\bar{b_e} + \beta_2 \bar{q_v}) \bm{\nabla} \cdot (\delta D \psi) \text{d}x
    \\
    + \frac{\alpha \Delta t}{2} \int_{\Gamma} \llbracket \delta D \bm{\psi}\rrbracket \langle \bar{b_e} + \beta_2 \bar{q_v} \rangle \text{dS}  = - \int_\Omega \bm{\psi} \cdot \bm{u}_r \text{d}x, \quad \forall \bm{\psi} \in \mathbb{V}^1,
    \\
    \int_\Omega \phi \delta D \text{d}x + \Delta t \int_\Omega \phi \bar{D} \bm{\nabla} \cdot \delta \bm{u} \text{d}x = - \int_\Omega \phi D_r \text{d}x, \quad \forall \phi \in \mathbb{V}^2,
    \\
    \int_{\Omega} \lambda \delta b_e \text{d}x -  \Delta t \int_{\Omega} \bar{b_{e}} \bm{\nabla} \cdot (\lambda {\delta \bm{u}}) \text{d}x + \Delta t \int_{\Gamma} \Tilde{\bar{b_e}} \llbracket \lambda {\delta \bm{u}}\rrbracket \text{d}S = - \int_\Omega \lambda b_{e_r} \text{d}x, \quad \forall \lambda \in \mathbb{V}^2,
\end{split}
\end{align}
where the increment we are solving for, $\delta \chi$, is given by $\delta \chi = (\delta \bm{u}, \delta D, \delta b_e)$. We choose the off-centring parameter $\alpha$ as 0.5. As in \citet{Melvin2024}, no off-centring parameter appears in the scalar equations. $\bar{q_v}$ is diagnosed using the saturation function at each stage from $q_t$ at the previous stage. We solve this system by eliminating the buoyancy equation via:
\begin{equation} \label{eq:b_e_expression}
    \delta b_e = -\Delta t (\delta \bm{u} \cdot \bm{\nabla}) \bar{b_e} + b_r,
\end{equation}
to give a reduced system for $\bm{u}$ and $D$. This reduced system is solved using a hybridisation method (described in detail in \citet{Firedrakebook} and \citet{Gibson2020slate, Bendall2020compatible}). The technique involves relaxing the continuity requirements on the normal components of functions in $\mathbb{V}^1$ and enforcing continuity through Lagrange multipliers as part of the solution formulation. We then use the solution $(\delta \bm{u}, \delta D)$ to reconstruct $\delta b_e$ via \eqref{eq:b_e_expression}.

\subsubsection{Moist linear solver}
\label{section:moist_linear_solver}
The moist linear operator $\mathcal S$ is derived by linearising the moist thermal shallow water equations about the reference state $\bm{0}, \bar{D}, \bar{b}, \bar{q_v}, \bar{q_c}$. The linearisation ignores topography and assumes that there is both water vapour and cloud present. We write the source term from the physics scheme as $P$, and linearise it by assuming that there is always some cloud to write:
\begin{equation}
    P = q_v - q_{sat}(q_v, D, b, q_v).
   \end{equation}
To linearise $P$ we derive the perturbation to $q_{sat}$ from a Taylor expansion of $q_{sat}(D, b, q_v)$ about $q_{sat}(\bar{D}, \bar{b}, \bar{q_v})$. This gives the expression for the linearised $P$ as:
\begin{equation}
    P = q_v + q_{sat}(\bar{D}, \bar{b}, \bar{q_v}) \left(\frac{1}{\bar{D}} D + \frac{\nu}{g} b - \frac{\beta_2 \nu}{g} q_v \right). \label{eq:P'_definition} 
\end{equation}
Then the weak form for the matrix-vector problem \eqref{eq:linear_solve} is given by:
\begin{align}
    \begin{split}
    & \int_\Omega {\bm{\psi}} \cdot \delta \bm{u} \text{d}x 
    + \alpha \Delta t \int_{\Omega} {\bm{\psi}} \cdot (f \bm{\hat{k}} \times \delta \bm{u}) \text{d}x
    - \alpha \Delta t \int_\Omega \bm{\nabla} \cdot (\bar{b} {\bm{\psi}}) \delta D \text{d}x
    + \alpha \Delta t \int_{\mathcal T} [[\bm{\psi}\bar{b}]] \langle \delta D \rangle \text {d}S
    \\
    &- \frac{\alpha \Delta t}{2}\int_\Omega \bar{D} \delta b \bm{\nabla} \cdot {\bm{\psi}} \text{d}x
    - \frac{\alpha \Delta t}{2}\int_\Omega \bar{b} \bm{\nabla} \cdot ( \delta D {\bm{\psi}}) \text{d}S
    + \frac{\alpha \Delta t}{2} \int_\mathcal{T} [[\delta D {\bm{\psi}}]] \langle \bar{b} \rangle \text{d}S
    = - \int_\Omega \bm{\psi} \cdot \bm{u}_r \text{d}x \quad \forall \bm{\psi} \in \mathbb{V}^1,
    \\
    &\int_\Omega \phi \delta D \text{d}x
    + \Delta t \int_\Omega \phi \bar{D} \bm{\nabla} \cdot \delta {\bm{u}} \text{d}x = - \int_\Omega \phi D_r \text{d}x \quad \forall \phi \in \mathbb{V}^2,
    \\
    &\int_\Omega \lambda \delta b \text{d}x - \Delta t \int_\Omega \bar{b} \bm{\nabla} \cdot (\lambda \delta \bm{u}) \text{d}x + \Delta t \int_\Gamma \langle \bar{b} \rangle \llbracket \lambda \delta \bm{u} \rrbracket \text{d}S
    + \Delta t \int_\Omega \lambda \beta_2 P \text{d}x = - \int_\Omega \lambda b_r \text{d}x \quad \forall \lambda \in \mathbb{V}^2,
    \\
    &\int_\Omega \tau_1 \delta q_v \text{d}x - \Delta t \int_\Omega \bar{q_v} \bm{\nabla} \cdot (\tau_1 \delta \bm{u}) \text{d}x + \Delta t \int_\Gamma \langle \bar{q_v} \rangle \llbracket \tau_1 \delta \bm{u} \rrbracket \text{d}S
    + \Delta t \int_\Omega \tau_1 P \text{d}x
    = - \int_\Omega \tau_1 q_{vr} \text{d}x \quad \forall \tau_1 \in \mathbb{V}^2,
    \\
    &\int_\Omega \tau_2 \delta q_c \text{d}x - \Delta t \int_\Omega \bar{q_c} \bm{\nabla} \cdot (\tau_2 \delta \bm{u}) \text{d}x + \Delta t \int_\Gamma \langle \bar{q_c} \rangle \llbracket \tau_2 \delta \bm{u} \rrbracket \text{d}S
    - \Delta t \int_\Omega \tau_2 P \text{d}x
    = - \int_\Omega \tau_2 q_{cr} \text{d}x \quad \forall \tau_2 \in \mathbb{V}^2,
    \end{split}
\end{align}
where the increment we are solving for, $\delta \chi$, is given by $\delta \chi = (\delta \bm{u}, \delta D, \delta b, \delta q_v, \delta q_c)$ and the residual $\chi_{r}$ is given by $\chi_{r} = (\bm{u}_r, D_r, b_r, q_{vr}, q_{vr})$. Like in the integrated physics formulation of the solver, no off-centring parameter appears in the scalar equations following \citet{Melvin2024}, and we chose the off-centring parameter in the $\bm{u}$ equation as 0.5.

Our strategy is to compute $P$ by first extracting the reference states and using $\bar{b}$ and $\bar{q_v}$ to compute $\bar{b_e}$. We then compute the saturation function using Equation \eqref{eq:md_sat_function} and these reference states, and use this and the current values of the variables to build the expression for $P$ from Equation \eqref{eq:P'_definition}. This is the $P$ that appears in the weak form at each timestep. The monolithic system is then solved with the conjugate gradient method. 

\subsection{Physics in the semi-implicit quasi-Newton timestep} \label{section:physics_in_the_SIQN_step}
The standard approach to physics in the SIQN timestep that is described above in Section \ref{section:SIQN_scheme} is to treat physics explicitly and call the physics scheme at the end of the timestep, after exiting both the inner and outer loops and just before the timestep is advanced. We are interested in experimenting with the approach to physics in the timestep and investigating whether or not there are any impacts from varying where in the timestep the physics scheme is called. In addition to the standard end-of-timestep physics which we will refer to as final physics, we have considered three other possibilities for physics:
\begin{enumerate}
    \item pre-loop physics: physics at the beginning of the timestep, before explicit forcing;
    \item outer loop physics: physics in the outer loop, after explicit transport; and
    \item inner loop physics: physics in the inner loop, treated semi-implicitly.
\end{enumerate}
The procedure for pre-loop physics is outlined in pseudocode in Algorithm \ref{algorithm:SIQN_slow_physics}. The physics increment is calculated directly after the first step of the algorithm, which as before is setting the guess for $\chi^{n+1}$ to $\chi^n$. Physics is then evaluated on the state at the current time level $\chi^n$ and returns $\chi_{pre\_loop}$, which is the state passed to the explicit forcing step. The forcing step acts on $\chi_{pre\_loop}$ to increment $\chi^n$, returning $\chi_{f,e}$. In the next step the algorithm enters the outer iterative loop and from there on it proceeds as before (see Section \ref{section:SIQN_scheme}). 

\begin{algorithm} 
	\caption{Pre-loop physics in the semi-implicit quasi-Newton time-stepping scheme} \label{algorithm:SIQN_slow_physics}
	\begin{algorithmic}[1]
		\State Set: $\chi^{n+1} = \chi^{n}$
        \State Pre-loop physics: $\chi_{pre\_loop} = \mathcal{P}(\chi^n)$
        \State Explicit forcing: $\chi_{f, e} = \chi_{pre\_loop} - \frac{1}{2}\Delta t \mathcal{F}(\chi^n)$
		\State Outer loop: \textbf{do}
			\State \hspace{\algorithmicindent} Update advecting velocity: $\bar{\bm{u}} = \frac{1}{2} (\bm{u}^{n+1} + \bm{u}^n)$
			\State \hspace{\algorithmicindent} Explicit transport: $\chi_T = \mathcal{A}_{\bar{\bm{u}}}(\chi_{f, e})$
            \State \hspace{\algorithmicindent} Inner loop: \textbf{do}
                \State \hspace{\algorithmicindent}\hspace{\algorithmicindent} Implicit forcing: $\chi_{f, i} = \chi_T - \frac{1}{2}\Delta t \mathcal{F}(\chi^{n+1})$
                \State \hspace{\algorithmicindent}\hspace{\algorithmicindent} Write residual: $\chi_{r} = \chi_{f, i} - \chi^{n+1}$
                \State \hspace{\algorithmicindent}\hspace{\algorithmicindent} Solve: $\mathcal{S}(\delta \chi) = \chi_{r}$ for $\delta \chi$
                \State \hspace{\algorithmicindent}\hspace{\algorithmicindent} Increment: $\chi^{n+1} = \chi^{n+1} + \delta \chi$
		      \State \hspace{\algorithmicindent} \textbf{End do}
        \State \textbf{End do}	
		\State Advance timestep: $\chi^n = \chi^{n+1}$
	\end{algorithmic} 
\end{algorithm}

The outer loop physics approach is described in Algorithm \ref{algorithm:SIQN_fast_physics}. The SIQN algorithm proceeds as usual (as described in Section \ref{section:SIQN_scheme}) until after explicit forcing in the outer loop. At that point physics is evaluated on the transported state $\chi_T$, returning $\chi_{outer\_loop}$. We then define the impact of this physics evaluation, $\chi_{rhs\_phys}$, as the difference between the transported state and the state after physics. This appears in the pseudocode as $\chi_{rhs\_phys} = \chi_{outer\_loop} - \chi_T$. The algorithm then enters the inner loop, where the first step of evaluating implicit forcing is the same as in the standard SIQN algorithm (so forcing acts on the transported state rather than the state with physics). As before we then define the residual $\chi_r$ as the difference between this newly-forced state and the best guess for the state at the next timestep, $\chi^{n+1}$. The difference here in comparison to the standard final physics approach is that at this stage we also add the impact of the physics evaluation to the residual, $\chi_{r} = \chi_{r} + \chi_{rhs\_phys}$. This means that the residual passed to the implicit solve in the next step includes the impact of physics. The solve method is applied as before to compute increments $\delta \chi$ which are added to the current guess to update $\chi^{n+1}$. After performing a fixed number of iterations of both the outer and inner loops we exit both loops before the timestep is advanced in the last step of the algorithm.

\begin{algorithm} 
	\caption{Outer loop physics in the semi-implicit quasi-Newton time-stepping scheme} \label{algorithm:SIQN_fast_physics}
	\begin{algorithmic}[1]
		\State Set: $\chi^{n+1} = \chi^{n}$
        \State Explicit forcing: $\chi_{f, e} = \chi^n - \frac{1}{2}\Delta t \mathcal{F}(\chi^n)$
		\State Outer loop: \textbf{do}
			\State \hspace{\algorithmicindent} Update advecting velocity: $\bar{\bm{u}} = \frac{1}{2} (\bm{u}^{n+1} + \bm{u}^n)$
			\State \hspace{\algorithmicindent} Explicit transport: $\chi_T = \mathcal{A}_{\bar{\bm{u}}}(\chi_{f, e})$
            \State \hspace{\algorithmicindent} Outer loop physics: $\chi_{outer\_loop} = \mathcal{P}(\chi_{T})$
            \State \hspace{\algorithmicindent} Physics increment: $\chi_{rhs\_phys} = \chi_{outer\_loop} - \chi_{T}$
            \State \hspace{\algorithmicindent} Inner loop: \textbf{do:}
                \State \hspace{\algorithmicindent}\hspace{\algorithmicindent} Implicit forcing: $\chi_{f, i} = \chi_T - \frac{1}{2}\Delta t \mathcal{F}(\chi^{n+1})$
                \State \hspace{\algorithmicindent}\hspace{\algorithmicindent} Write residual: $\chi_{r} = \chi_{f, i} - \chi^{n+1}$
                \State \hspace{\algorithmicindent}\hspace{\algorithmicindent} Add physics increment to residual: $\chi_{r} = \chi_{r} + \chi_{rhs\_phys}$ 
                \State \hspace{\algorithmicindent}\hspace{\algorithmicindent} Solve: $\mathcal{S}(\delta \chi) = \chi_{r}$ for $\delta \chi$
                \State \hspace{\algorithmicindent}\hspace{\algorithmicindent} Increment: $\chi^{n+1} = \chi^{n+1} + \delta \chi$
			\State \hspace{\algorithmicindent} \textbf{End do}
        \State \textbf{End do}
		\State Advance timestep: $\chi^n = \chi^{n+1}$
	\end{algorithmic} 
\end{algorithm}

Inner loop physics treats physics in the same way as forcing; that is, physics is evaluated both explicitly in the outer loop and implicitly in the inner loop. The approach can be written as:
\begin{equation}
    \chi^{n+1} - \frac{1}{2} \Delta t \mathcal{F}(\chi^{n+1}) + \beta \Delta t \mathcal{P} (\chi^{n+1}) = \mathcal{A}_{\bar{\bm{u}}} \left[ \chi^n + \frac{1}{2} \Delta t \mathcal{F}(\chi^{n}) - (1 - \beta)\Delta t \mathcal{P} (\chi^{n}) \right],
\end{equation}
where $\beta$ is the semi-implicit off-centring parameter for the physics. $\beta$ takes a value between 0 and 1, where $\beta=1$ corresponds to a fully implicit evaluation and $\beta=0$ zero a fully explicit evaluation. The linear system that is solved in the inner loop with this approach is the moist system described in Section \ref{section:moist_linear_solver}, as opposed to the dry linear solver which is used with all other approaches to physics. The algorithm for inner loop physics is given in Algorithm \ref{algorithm:SIQN_inner_physics}.

\begin{algorithm} 
	\caption{Inner loop physics in the semi-implicit quasi-Newton time-stepping scheme} \label{algorithm:SIQN_inner_physics}
	\begin{algorithmic}[1]
		\State Set: $\chi^{n+1} = \chi^{n}$
        \State Explicit forcing: $\chi_{f, e} = \chi^n - \frac{1}{2} \Delta t \mathcal{F}(\chi^n)$
        \State Explicit physics increment: $\chi_{exp\_incr} = \mathcal{P}(\chi^n) - \chi^n$
        \State Add explicit physics increment: $\chi_{p,e} = \chi_{f,e} + (1 - \beta)\Delta t \chi_{exp\_incr}$
		\State Outer loop: \textbf{do}
			\State \hspace{\algorithmicindent} Update advecting velocity: $\bar{\bm{u}} = \frac{1}{2} (\bm{u}^{n+1} + \bm{u}^n)$
			\State \hspace{\algorithmicindent} Explicit transport: $\chi_T = \mathcal{A}_{\bar{\bm{u}}}(\chi_{p, e})$
            \State \hspace{\algorithmicindent} Inner loop: \textbf{do}
                \State \hspace{\algorithmicindent}\hspace{\algorithmicindent} Implicit forcing: $\chi_{f, i} = \chi_T - \frac{1}{2} \Delta t \mathcal{F}(\chi^{n+1})$
                \State \hspace{\algorithmicindent}\hspace{\algorithmicindent} Write residual: $\chi_{r} = \chi_{f, i} - \chi^{n+1}$
                \State \hspace{\algorithmicindent}\hspace{\algorithmicindent} Implicit physics increment: $\chi_{imp\_incr} = \mathcal{P}(\chi^{n+1}) - \chi^{n+1}$
                \State \hspace{\algorithmicindent}\hspace{\algorithmicindent} Add implicit physics increment to residual: $\chi_{r} = \chi_{r} + \beta \Delta t \chi_{imp\_incr}$ 
                \State \hspace{\algorithmicindent}\hspace{\algorithmicindent} Solve: $\mathcal{S}(\delta \chi) = \chi_{r}$ for $\delta \chi$
                \State \hspace{\algorithmicindent}\hspace{\algorithmicindent} Increment: $\chi^{n+1} = \chi^{n+1} + \delta \chi$
			\State \hspace{\algorithmicindent} \textbf{End do}
        \State \textbf{End do}
		\State Advance timestep: $\chi^n = \chi^{n+1}$
	\end{algorithmic} 
\end{algorithm}

\section{Results}
In this section we demonstrate both the new integrated-physics formulation of the moist thermal shallow water equations and the new moist linear solver, and we give results from physics-dynamics coupling experiments in the SIQN timestepper. We begin by verifying that the new integrated physics formulation gives the expected results using a steady-state geostrophically balanced test. We then use a moist gravity wave test case to demonstrate the new formulation's capabilities as a `ground-truth' model to compare coupling approaches in the split-physics moist thermal formulation against. Both of these tests involve a careful initialisation procedure to ensure that the test is set up in the same way in both formulations, so that comparing results between the integrated-physics and the split-physics moist thermal shallow water models is `fair'. After demonstrating the suitability of the model as a `ground truth' we will use it as a reference solution to compare approaches to physics-dynamics coupling in the SIQN timestepper.

All of our results are produced using Gusto, the dynamical core toolkit built on the Firedrake finite element library. Firedrake provides automated code generation for the solution of partial differential equations using the finite element method \citep{FiredrakeUserManual}. Code is automatically generated using the Unified Form Language \citep{UFL} and equations are provided to PETSc \citep{PETSc}, which provides direct access to runtime configurable iterative solvers and preconditioners. In this work all experiments use an icosahedral grid. 

\subsection{Setting up a fair test} \label{section:iterative_procedure_for_md_formulation}
The essence of a fair comparison between the integrated-physics formulation and the split-physics formulation is that both formulations begin with the same initial conditions and that moist physics processes are triggered by the same conditions in both formulations, even though the buoyancy and moisture prognostic fields vary between the two formulations. Our technique is to begin by specifying an initial value for the buoyancy $b$ and vapour $q_v$, and to use these to compute the initial equivalent buoyancy for the integrated-physics formulation via:
\begin{equation} \label{eq:b_be_relationship}
    b_e = b - \beta_2 q_v.
\end{equation}
Frequently, the initial vapour will depend on the initial saturation function as we wish to begin a test with an amount of vapour that is some proportion of the saturation function; for example, at saturation or just below it. This means we need to solve
\begin{equation}
    b_e = b - \beta_2 q_{sat},
\end{equation}
where $q_{sat}$ depends on $b_e$. We solve this non-linear problem using a Newton-Raphson method. Further details of this initialisation procedure are given in Appendix \ref{app:set_up}.

Then to set up a test in both formulations we specify an initial velocity and depth that are the same in both moist shallow water formulations. We define an initial guess for the buoyancy $b$ and the vapour $q_v$. We then use the iterative procedure outlined in Appendix \ref{app:set_up} to compute the initial saturation function $q_{sat}$ and thus the initial vapour, which in turn allows us to define the initial equivalent buoyancy $b_e$. In the integrated-physics formulation we prescribe the initial total moisture $q_t$ rather than the vapour $q_v$, but typically $q_t$ is based on $q_v$. If the test is in the split-physics formulation we back out the initial buoyancy $b$ using the $b_e$ and $q_v$ computed from the iterative procedure via Equation \eqref{eq:b_be_relationship}.

This initialisation procedure means that tests in both the integrated-physics and the split-physics formulations of the moist thermal shallow water equations have equivalent initial conditions. To ensure that moist processes are triggered by the same conditions in both formulations we compute the saturation function in both cases based on the equivalent buoyancy $b_e$. In the split-physics formulation this means converting the buoyancy to equivalent buoyancy using the relationship in Equation \eqref{eq:b_be_relationship} before passing this equivalent buoyancy as an argument to the saturation function. We also choose the $\gamma_v$ parameter in the split-physics moist physics scheme to be 1 so that the full amount of the moist species is converted in a phase change, rather than the proportion given by Equation \eqref{eq:gamma_v}. This means that as well as matching the initial conditions in both formulations and triggering a phase change with the same conditions in both formulations, a phase change also produces the same output in both formulations.

\subsection{Model Validation}

We begin by validating that the new integrated physics formulation gives the expected convergence results and that the split-physics formulation of the equations converges to the integrated physics formulation as the timestep size is reduced. These experiments all use the dry linear solver (which excludes moisture from the linear solve stage of the SIQN scheme). 

\subsubsection{Over-saturated steady state geostrophic test} \label{section:oversaturated_W2}
The first test of the new integrated-physics formulation is a steady state geostrophic test case, which verifies that the model can maintain a steady state and builds confidence in our formulation. We run the test using both the split-physics and the integrated-physics formulations of the moist thermal shallow water equations. The experiment is a version of the geostrophically-balanced test from \citet{Hartney2025} which is based on the second test of the \citet{Williamson1992} test suite and test 2 of \citet{ZA2015}. Unlike the approach in \citep{Hartney2025}, however, in this test we begin with an over-saturated state so that we are starting with vapour at saturation as well as some cloud. Because the test is a steady state, the results are assessed on how well the final state compares to the initial state and we produce convergence plots of the difference between the final and initial states as the resolution is increased.  

We take rotation rate of the Earth $\Omega = 7.292 \times 10^{-5} \text{s}^{-1}$, acceleration due to gravity $g = 9.80616 \text{ m } \text{s}^{-2}$ and radius of the Earth $R = 6371220 \text{ m}$. The initial velocity and depth are the same in both formulations of this test and are given by (with $(\lambda, \phi)$ representing spherical (longitude, latitude) coordinates):
\begin{align}
    &\ u(\lambda, \phi) = u_0\cos\phi, \\
    &\ v(\lambda, \phi) = 0, \\
    &D(\lambda, \phi) = H - \frac{1}{g} (\omega + \sigma) \text{sin}^2\phi,
\end{align}
with $u$ and $v$ the zonal and meridional components of $\bm{u}$ respectively, $u_0 = 20 \text{ m s}^{-1}$, $\omega = \left( \Omega R u_0 + \frac{u_0^2}{2}\right)$, $\sigma = \omega/10$, and the background depth $H$ is $\Phi_0/g$ where $\Phi_0 = 3 \times 10^4 \text{ m s}^{-2}$ and $g$ is the acceleration due to gravity. The guess for the buoyancy that is used in the iterative procedure to find $b_e$ and the saturation function is the initial buoyancy from the steady state test in Section 5.1 of \citep{Hartney2025}:
\begin{equation}
    b_0(\lambda, \phi) = g \left( 1 - \frac{\theta_0 + \sigma \text{cos}^2\phi \left[ (\omega + \sigma) \text{cos}^2\phi + 2(\Phi_0 - \omega - \sigma)\right]}{\Phi_0^2 + (\omega + \sigma)^2 \text{sin}^4 \phi - 2 \Phi_0 (\omega + \sigma) \text{sin}^2 \phi}\right),
    \label{eq:md_W2_buoyancy_guess}
\end{equation}
where $\theta_0 = \epsilon \Phi_0^2$ with $\epsilon = 1/300$. The initial guess for $q_v$ for the iterative procedure is taken as 0.02. Then the initial $b_e$ comes from the outcome of the iterative procedure, and both the initial vapour and cloud are set equal to the initial saturation function, which is also obtained from the iterative procedure. The parameters for the saturation function are $\nu = 1.5$ and $q_0 = 0.007$.

\begin{figure}[htp]
    \centering
    \includegraphics[width=0.75\textwidth]{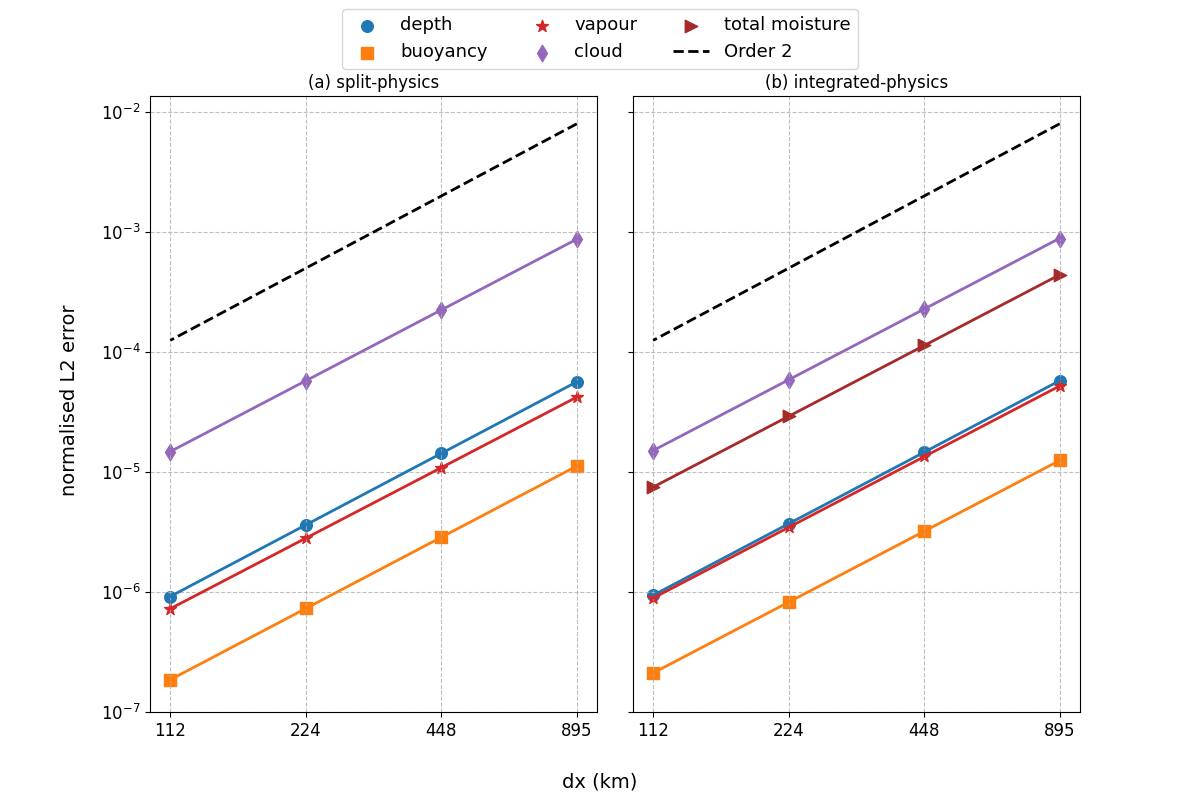}
    \caption[Convergence in the split-physics and integrated-physics oversaturated Williamson 2 test]{$L^2$ error of the over-saturated moist Williamson 2 test in the split-physics and integrated-physics frameworks as a function of spatial resolution $dx$. (a) the split-physics formulation, and (b) the integrated-physics formulation. The $L^2$ norm of the error is normalised by dividing it by the $L^2$ norm of the field in each case. Second-order convergence lines are also shown on the plot for comparison.}
    \label{fig:oversaturated_W2_convergence}
\end{figure}

Figure \ref{fig:oversaturated_W2_convergence} shows convergence plots of the $L^2$ error norm of each field after 5 days in both the split-physics formulation and the integrated-physics formulation. The timestep is varied as the spatial resolution changes so as to maintain a constant advective Courant number of approximately 0.1. Because of the polynomial degree used in the finite elements we expect second-order convergence from our time-stepping scheme. Second-order convergence lines are shown on all plots, showing that all fields in both formulations converge at the expected rate. The results demonstrate that the integrated-physics formulation gives the expected behaviour for this steady-state test. 

\subsubsection{Moist gravity wave test} \label{section:moist_gravity_wave_test}
The second test verifies that the split-physics formulation of the moist thermal shallow water model will converge to the integrated-physics formulation as the timestep is reduced. The aim of this test is to build the case for the model's suitability as a reference solution with which to compare coupling strategies in the split-physics moist thermal model against. The test should challenge the model more than the steady state test, but should also be sufficiently smooth so that we can expect convergence behaviour as the timestep is reduced. We have devised a shallow water version of a moist gravity wave test which fulfils both of these criteria. The fact that the test does not have an analytic solution is immaterial for our purpose, which is simply to show that the split-physics solution converges to the integrated-physics solution in the limit of decreasing timestep.     

The test begins with a balanced flow to which a perturbation to the depth field is added. The perturbation triggers a propagating gravity wave. In both the integrated-physics and the split-physics version of the test we begin with a uniform amount of total moisture, which is partitioned between vapour and cloud by the saturation function. The initial balanced velocity and depth fields are the same as in the steady state test of the previous section. The perturbation to the depth field has the same form as the mountain in the flow over a mountain test from the \citet{Williamson1992} test suite:
\begin{equation}
    p = h_0\left(1 - \frac{1}{R} \text{min} \left[ R, \sqrt{(\lambda - \lambda_c)^2 + (\phi - \phi_c)^2} \right]\right),
\end{equation}
with $h_0 = 2000$ m, $R = \pi/9$, $\lambda_c = 3\pi/2$ and $\phi_c = \pi/6$.
This perturbation is added to the initial depth field to trigger the gravity waves. The initial velocity (in terms of zonal and meridional components $(u,v)$) and depth fields are given by:
\begin{align}
    &u(\lambda, \phi) = u_0 \text{cos} \phi, \\
    &v(\lambda, \phi) = 0, \\
    &D(\lambda, \phi) = H - \frac{1}{g} (\omega + \sigma) \text{sin}^2\phi + p,
\end{align}
where, as in the steady state test, $u_0 = 20 \text{ m s}^{-1}$, $\omega = \left( \Omega R u_0 + \frac{u_0^2}{2}\right)$, $\sigma = \omega/10$, and the background depth $H$ is $\Phi_0/g$ where $\Phi_0 = 3 \times 10^4 \text{ m s}^{-2}$ and $g$ is the acceleration due to gravity. The guess for the buoyancy that is used in the iterative procedure to find $b_e$ and the saturation function is the same as that used in the steady state test. The initial guess for $q_v$ for the iterative procedure is taken as 0.03 and the parameters for the saturation function are $\nu = 1.5$ and $q_0 = 0.0115$. The initial total moisture in the integrated-physics formulation $q_t$ is 0.03. The initial vapour $q_v$ in the split-physics formulation is set equal to the saturation function $q_{sat}$, and the initial cloud field $q_c$ is the difference between 0.03 and the saturation function. 
We run the test for 5 days.

We take a high-resolution version of the test in the integrated-physics formulation to act as a reference solution, and then vary the timestep in the split-physics version of the test to assess how well the split-physics solution converges to the reference solution as the timestep is reduced. All tests use the same icosahedral grid that corresponds to a spatial resolution of approximately 224 km. The high resolution reference uses a timestep of 50 seconds and the split-physics experiments use timesteps varying from 800 seconds to 100 seconds. To compare solutions between formulations we convert the equivalent buoyancy in the integrated-physics reference solution to buoyancy via:
\begin{equation}
    b = b_e + \beta_2 q_v. \label{eq:equivalent_b_relationship}
\end{equation}
We can then compare fields in the split-physics formulation to the integrated-physics formulation using the $L^2$ norm of the difference.

\begin{figure}
    \centering
    \includegraphics[width=\linewidth]{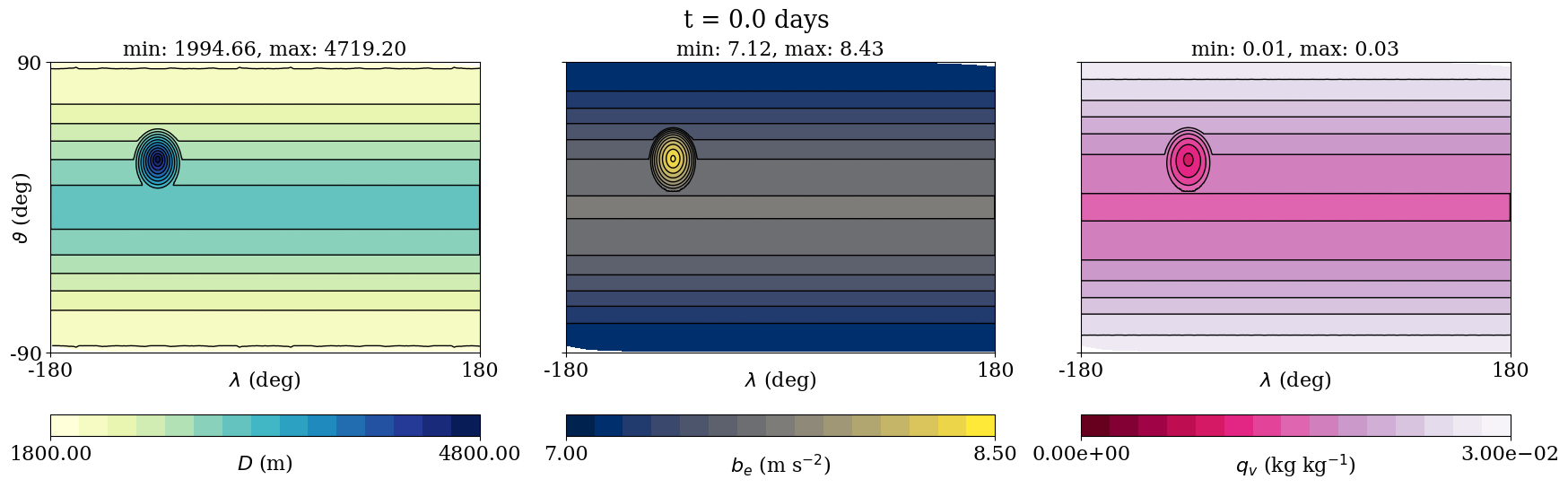}
    \caption{Initial conditions for the depth ($D$), equivalent buoyancy ($b_e$) and vapour ($q_v$) in the integrated-physics formulation of the moist thermal gravity wave test.}
    \label{fig:gw_ref_ICs}
\end{figure}

\begin{figure}
    \centering
    \includegraphics[width=\linewidth]{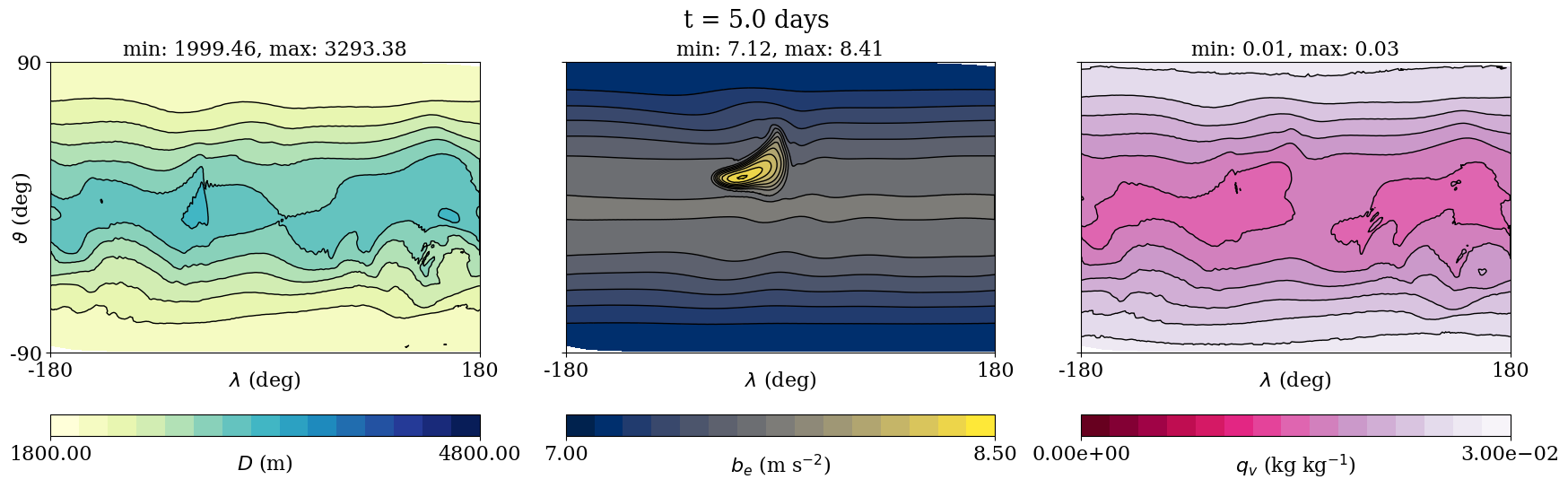}
    \caption{Depth ($D$), equivalent buoyancy ($b_e$) and vapour ($q_v$) after 5 days in the integrated-physics formulation of the moist thermal gravity wave test.}
    \label{fig:gw_ref_final}
\end{figure}

\begin{figure}
    \centering
    \includegraphics[width=\linewidth]{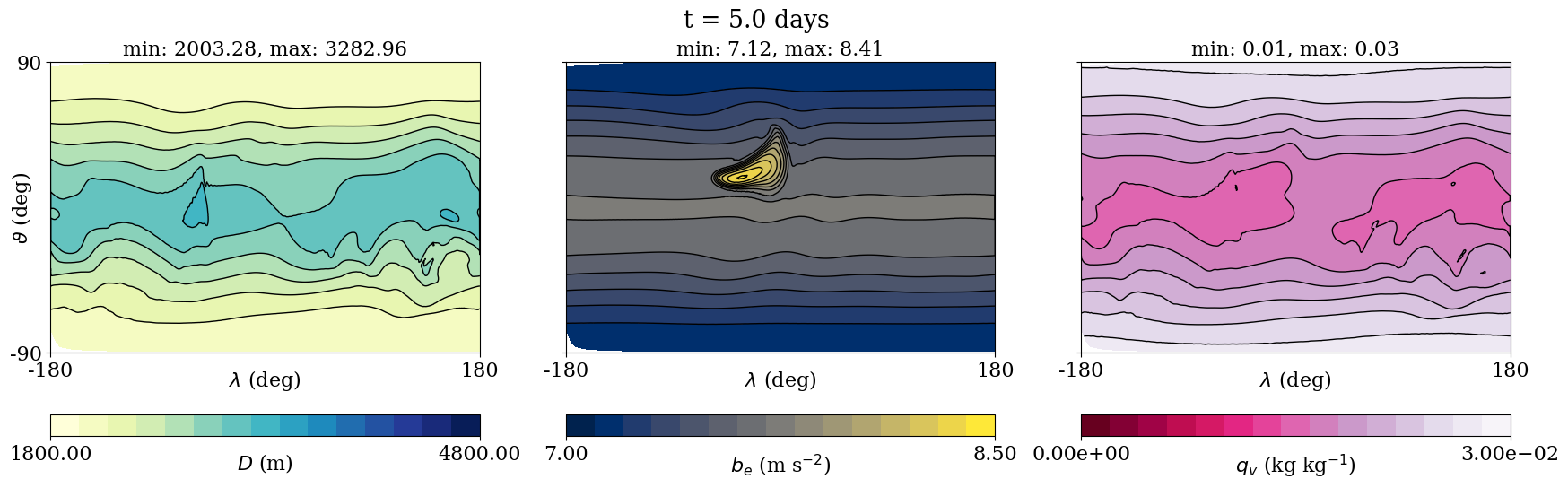}
    \caption{Depth ($D$), equivalent buoyancy ($b_e$) and vapour ($q_v$) after 5 days in the split-physics formulation of the moist thermal gravity wave test. The equivalent buoyancy comes from converting the buoyancy in the test to equivalent buoyancy via Equation \eqref{eq:equivalent_b_relationship}.}
    \label{fig:gw_final_phys_final}
\end{figure}

Figure \ref{fig:gw_ref_ICs} shows the set up for the test with the initial depth, equivalent buoyancy and vapour from the integrated-physics formulation, and Figures \ref{fig:gw_ref_final} and \ref{fig:gw_final_phys_final} show the same fields after they have evolved for 5 days in the integrated-physics formulation and the split-physics formulation, respectively.

\begin{figure}[htp]
    \centering
    \includegraphics[width=0.65\textwidth]{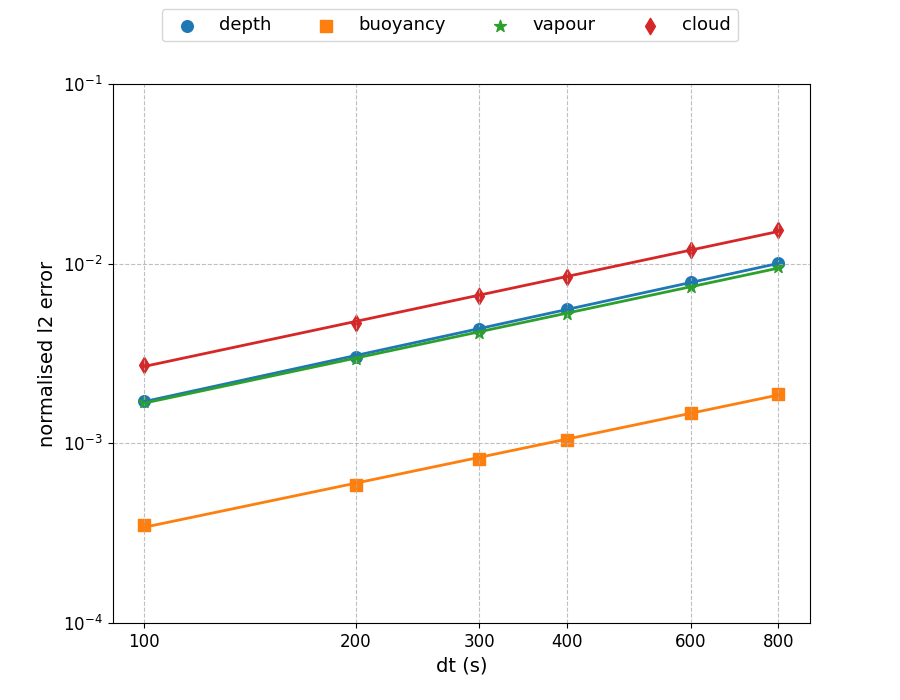}
    \caption[Convergence of the split-physics formulation of the moist thermal gravity wave test to the integrated-physics solution]{The $L^2$ error of the split-physics moist thermal solution for the moist gravity wave test with varying timesteps ($dt$) in seconds, as compared to the integrated-physics reference solution (with a timestep of 50 s). The error in each field is normalised by dividing it by the norm of the reference solution field.}
    \label{fig:moist_gravity_wave_convergence}
\end{figure}

Figure \ref{fig:moist_gravity_wave_convergence} shows the normalised $L^2$ error in the split-physics formulation of the moist thermal gravity wave test, where the integrated-physics solution is treated as a reference solution. The errors are normalised by dividing them by the norm of the reference solution for each field. The vapour and cloud in the integrated-physics formulation are computed by using the saturation function to partition the total moisture. The error in all fields decreases as the timestep is reduced, demonstrating how the split-physics solution converges to the integrated-physics one in the limit of infinitely-small timestep size, and hence demonstrating the integrated-physics model's suitability as a reference solution.  

\subsection{Physics in the semi-implicit quasi-Newton timestepper}
\begin{figure}
    \centering
    \includegraphics[width=\linewidth]{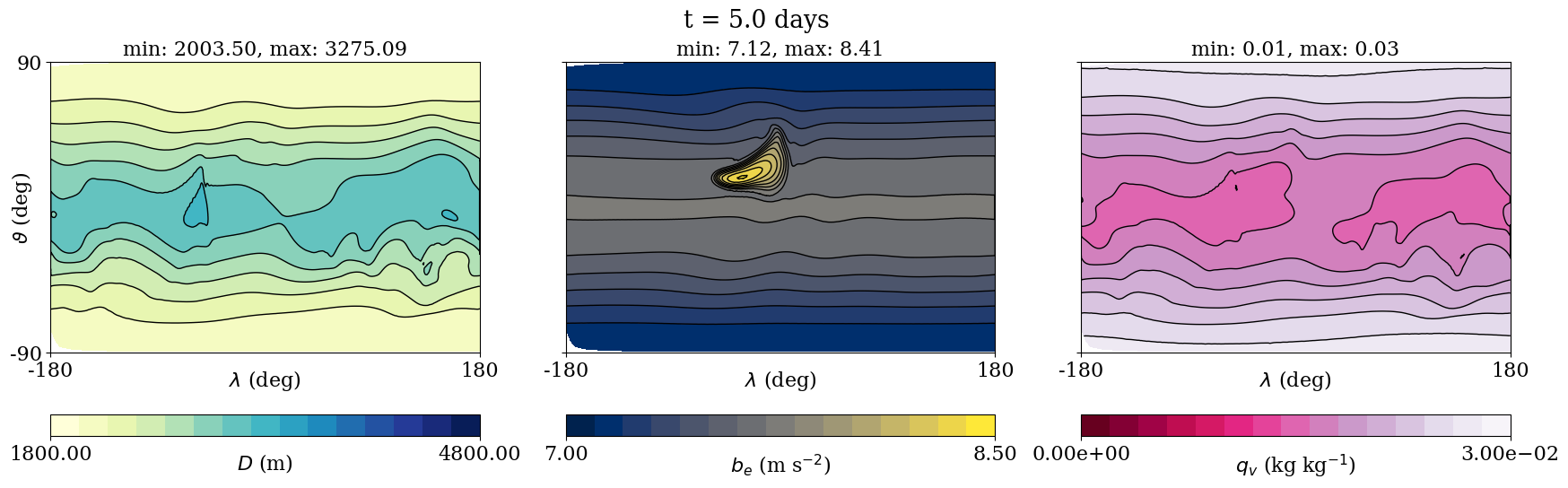}
    \caption{Depth ($D$), equivalent buoyancy ($b_e$) and vapour ($q_v$) after 5 days in the split-physics formulation of the moist thermal gravity wave test, where physics is evaluated with the pre-loop physics approach. The equivalent buoyancy comes from converting the buoyancy in the test to equivalent buoyancy via Equation \eqref{eq:equivalent_b_relationship}.}
    \label{fig:gw_slow_phys_final}
\end{figure}

\begin{figure}
    \centering
    \includegraphics[width=\linewidth]{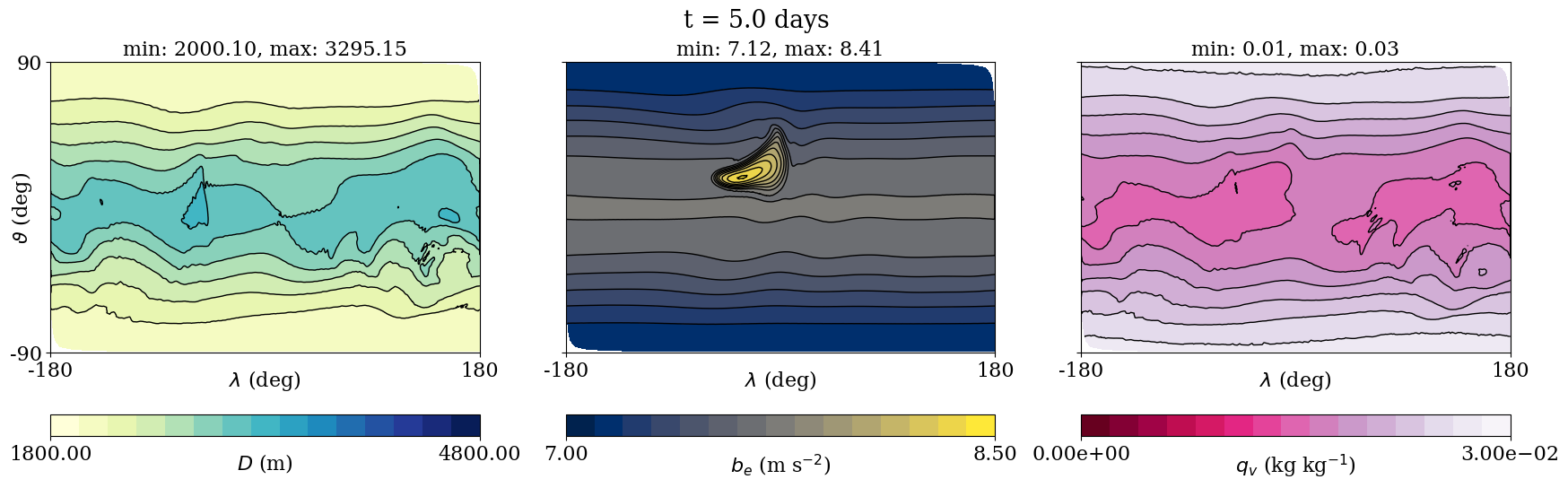}
    \caption{Same as Figure \ref{fig:gw_slow_phys_final}, but where physics is evaluated in the outer loop.}
    \label{fig:gw_fast_phys_final}
\end{figure}

\begin{figure}
    \centering
    \includegraphics[width=\linewidth]{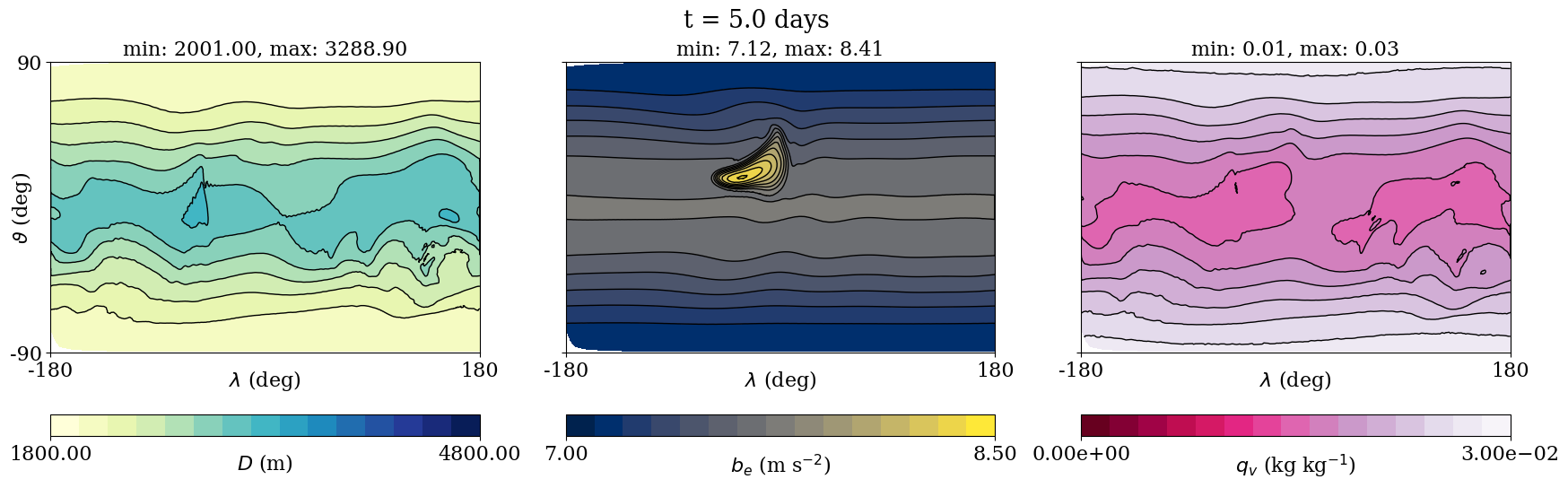}
    \caption{Same as Figure \ref{fig:gw_slow_phys_final}, but where physics is evaluated in the inner loop.}
    \label{fig:gw_inner_phys_beta1_final}
\end{figure}

In this experiment we investigate different approaches to physics-dynamics coupling in the SIQN time-stepping scheme. As described in Section \ref{section:SIQN_scheme}, the approach to physics in the split-physics formulation we have taken thus far has been to treat physics explicitly and to call it at the end of the timestep (so-called final physics). Here we use the other approaches to physics in the SIQN timestep described in Section \ref{section:physics_in_the_SIQN_step}, and test them using the moist thermal gravity wave test described in the previous section. The fields after 5 days for pre-loop physics, outer loop physics and inner loop physics (with the semi-implicit off-centring parameter $\beta = 1.0$) are given in Figures \ref{fig:gw_slow_phys_final}, \ref{fig:gw_fast_phys_final} and \ref{fig:gw_inner_phys_beta1_final}, respectively.


We next assess how the error in the solutions vary between physics approaches, by examining the differences between the split-physics solutions and the integrated-physics reference solution. We plot the normalised $L^2$ error of the fields after 5 days using the integrated-physics formulation as a reference solution, with the same constant spatial resolution as the experiment in Section \ref{section:oversaturated_W2} (224 km cell lengths on average) and varying the timestep. The integrated-physics reference solution uses the same resolution and a timestep $\Delta t$ of 50 s. We test inner loop physics both with the semi-implicit off-centring parameter $\beta = 1.0$ and $\beta = 0.5$.

\begin{figure}
    \centering
    \includegraphics[width=\linewidth]{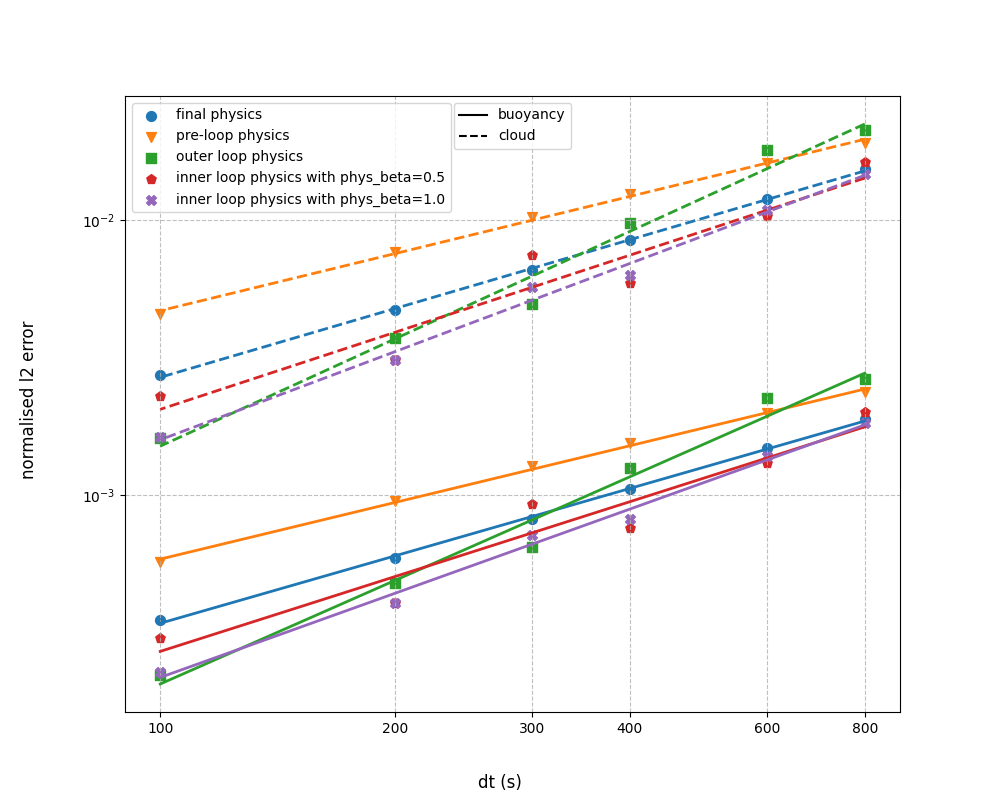}
    \caption[Comparison of physics approaches (by approach) in the moist gravity wave test]{Normalised $L^2$ error (as compared to an integrated-physics reference solution) for the moist gravity wave test using five different approaches to physics in the SIQN timestep. (a) physics at the end of the timestep, (b) physics at the beginning of the timestep (pre-loop physics), (c) physics in the outer loop, (d) physics in the inner loop (fully implicit), and (e) physics in the inner loop (semi-implicit).}
    \label{fig:gw_compare_physics}
\end{figure}

Figure \ref{fig:gw_compare_physics} shows the convergence rate of the error in both the buoyancy field and the cloud field for the different physics approaches. There is a clear distinction in the rate of convergence between the five different approaches. Outer loop physics converges the most quickly, but the overall error is generally lowest in the inner loop physics approaches. These are the approaches where the physics is most tightly coupled to the dynamics and thus we expect the lowest error from these schemes.

Another important observation from these results is the difference between final physics and pre-loop physics. Final physics converges slightly more quickly than pre-loop physics, but the more significant difference between the approaches is the size of the error, which is substantially bigger in pre-loop physics. Both final physics and pre-loop physics evaluate physics once per timestep outside both the inner and outer loop, and so could be expected to show similar errors and convergence rates. The difference can be explained by the fact that explicit forcing is evaluated on different states in these two approaches: the timelevel $n$ state $\chi^n$ in final physics, as compared to the state after physics $\chi_{pre\_loop}$ in pre-loop physics. This result could have interesting implications for physics-dynamics coupling strategies. For example, if the physics scheme is to be called only once per timestep (perhaps for operational or efficiency reasons), it appears that the coupling error is smaller when this call is made at the end of the timestep after increments from forcing have been computed rather than the beginning of the timestep. 

\section{Conclusions}
In this paper we have used two variations of the moist thermal shallow water equations to present a simplified framework that can be used to investigate questions around how physics parametrisations interact with the dynamical core in a weather or climate model. The shallow water equations are computationally cheap but still retain many relevant atmospheric effects. The inclusion of moisture in the system means a model with a coupled physics process, offering a simple and inexpensive tool to experiment with physics-dynamics coupling approaches.

The moist thermal shallow water equations include the effect of moist physics in the shallow water model by coupling the thermal shallow water equations to a moist physics scheme. In this work we introduced the integrated-physics moist thermal shallow water equations; a variation of the usual split-physics formulation of the moist shallow water equations where all the moist physics processes are captured in the equation dynamics. This removes the need for a physics-dynamics coupling approach, thus making the integrated-physics model suitable as a reference solution to compare different coupling strategies in the split-physics model to. We developed a test case based on a moist gravity wave that we ran in both formulations, and, using the integrated-physics solution as a reference, assessed the impact of different coupling approaches on the accuracy of the split-physics solution.

Following the approach of the next-generation Met Office dynamical core, GungHo, our models are discretised in time using the SIQN timestepping scheme. We use four different strategies to evaluate physics in the SIQN timestep: at the beginning of the timestep, at the end of the timestep, in the outer loop, and in the inner loop. To evaluate physics in the inner loop it was necessary to develop a version of the linear solver that includes moist prognostic fields, which is presented for the first time here. The results demonstrate clear differences between the coupling approaches, both in terms of the accuracy of the solution and the rate at which the solution converges to the reference. The results also demonstrate more broadly the usefulness of the moist shallow water model in general - and the integrated-physics formulation in particular - to advance understanding of physics-dynamics coupling issues.  

\newpage

\section*{Funding Acknowledgements}
N.H was supported by a Natural Environment Research Council (NERC) GW4+ Doctoral Training Partnership 2 Studentship (NE/S007504/1). For the purpose of open access, the authors have applied a creative commons attribution (CC BY) licence to any author accepted manuscript version arising from this submission.

\section*{Data Availability}
The data that support the findings of this study are available from the corresponding author upon reasonable request.

\section*{Conflicts of Interest}
The authors have no conflict of interest to declare.

\newpage
\appendix
\section{Initialisation procedure} \label{app:set_up}
Here we give details of the procedure used to initialise test cases in both the split-physics and integrated-physics formulation of the moist shallow water equations such that the initial conditions are the same in both formulations. 
This is not altogether straightforward when the buoyancy and moisture prognostic fields vary between the two formulations. Our technique is to begin by specifying an initial value for the buoyancy $b$ and vapour $q_v$, and to use these to compute the initial equivalent buoyancy for the integrated-physics formulation via Equation \eqref{eq:b_be_relationship}. Because the initial vapour will often depend on the initial saturation function, we have that $q_v$ will be some function of $q_{sat}$, which is a function of $b_e$ via Equation \eqref{eq:md_sat_function}. We have prescribed only the initial $b$ field, however, and to recover the $b_e$ field from the $b$ field we need to know $q_v$, which is the quantity we are trying to diagnose with the saturation function. In summary, the problem to solve for the initial value of $q_v$ is, with $\xi$ a parameter which dictates how close the initial vapour is to saturation: 
\begin{equation}
    q_v = (1 - \xi)q_{sat} = (1 - \xi )\frac{q_0 H}{D + B} e^{\nu(1-b/g+\beta_2 q_v/g)},
\end{equation}
and because $q_v$ appears on the right-hand side, we will solve this using an iterative procedure.

Our approach is to take $\xi = 0$ and to set:
\begin{equation}
    q_v = q_{sat}(b - \beta_2 q_v)
\end{equation}
and define $f(q_v) = q_{sat} - q_v$. We then solve $f(q_v) = 0$ for $q_v$ using a Newton-Raphson method, defined by:
\begin{equation}
    q_v^{n+1} = q_v^n - \frac{f(q_v^n)}{f'(q_v^n)},
\end{equation}
where $n$ is the iteration count and $f'$ represents the derivative of $f$ with respect to $q_v$:
\begin{align}
    f'(q_v) &= \frac{\partial}{\partial q_v} (q_{sat}) - \frac{\partial}{\partial q_v} (q_v) \\
    &= \frac{\nu \beta_2}{g} q_{sat} - 1
\end{align}
The Newton method then becomes:
\begin{equation}
    q_v^{n+1} = q_v^n - \frac{q_{sat}^n - q_v^n}{\frac{\nu \beta_2}{g} q_{sat}^n - 1},
\end{equation}
where $q_{sat}$ is updated at each iteration via:
\begin{equation}
    q_{sat}^n = \frac{q_0 H}{D + B} e^{\nu(1-b/g+\beta_2 q_v^n/g)}.
\end{equation}
This iterative procedure thus gives the initial saturation function (from which the initial vapour is computed) as well as the initial equivalent buoyancy. In practice we have found that 10 iterations are sufficient to converge the method, provided the initial guess for $q_v$ is good enough. 

\newpage
\bibliographystyle{unsrtnat}
\bibliography{paper2}

\end{document}